\documentclass[a4paper,12pt]{article}

\usepackage{hyperref}
\usepackage[pdftex]{graphicx}
\usepackage[section] {placeins}
\usepackage{enumerate}
\usepackage[normalem]{ulem}
\usepackage{amsmath}
\usepackage{amsthm}
\usepackage{latexsym}
\usepackage{amsfonts}
\usepackage{amssymb}
\usepackage{mathtools}
\usepackage{longtable}

\usepackage{tikz-cd}
\usetikzlibrary{graphs}
\usetikzlibrary{arrows.meta}
\usetikzlibrary{shapes.geometric}

\usepackage{authblk}

\usepackage{latexsym,amssymb,amsfonts, amsthm, amsmath}
\usepackage[english]{babel}
\usepackage{bm}
\usepackage{mathrsfs}
\usepackage{algpseudocode}
\usepackage{algorithm}

\algnewcommand{\IIf}[1]{\State\algorithmicif\ #1\ \algorithmicthen}
\algnewcommand{\EndIIf}{\unskip\ \algorithmicend\ \algorithmicif}

\newtheorem{thm}{Theorem}[section]
\newtheorem{lem}[thm]{Lemma}

\newtheorem{prop}[thm]{Proposition}

\newtheorem{rem}[thm]{Remark}
\newtheorem*{thm*}{Theorem \ref{th:main3}}
\newtheorem*{thm*2}{Theorem \ref{th:main1}}

\theoremstyle{definition}
\newtheorem{defn}[thm]{Definition}

\usepackage{fullpage}

\makeatletter
\newcommand{\subjclass}[2][1991]{%
\let\@oldtitle\@title%
\gdef\@title{\@oldtitle\footnotetext{#1 \emph{Mathematics subject classification.} #2}}%
}
\newcommand{\keywords}[1]{%
\let\@@oldtitle\@title%
\gdef\@title{\@@oldtitle\footnotetext{\emph{Key words and phrases.} #1.}}%
}
\makeatother

\begin{document}
\title{A Littlewood-Offord kind of problem in $\mathbb{Z}_p$\\ and $\Gamma$-sequenceability}
\author{Simone Costa}

\affil{DICATAM, Sez.~Matematica, Universit\`a degli Studi di Brescia, Via Branze~43, I~25123 Brescia, Italy}

\maketitle
\begin{abstract}
The Littlewood-Offord problem is a classical question in probability theory and discrete mathematics, proposed, firstly by Littlewood and Offord in the 1940s. Given a set $A$ of integer, this problem asks for an upper bound on the probability that a randomly chosen subset $X$ of $A$ sums to an integer $x$.

This article proposes a variation of the problem, considering a subset $A$ of a cyclic group of prime order and examining subsets $X\subseteq A$ of a given cardinality $\ell$. The main focus of this paper is then on bounding the probability distribution of the sum $Y$ of $\ell$ i.i.d. $Y_1,\dots, Y_{\ell}$ whose support is contained in $\mathbb{Z}_p$. The main result here presented is that, if the probability distributions of the variables $Y_i$ are bounded by $\lambda \leq 9/10$, then, assuming that $p> \frac{2}{\lambda}\left(\frac{\ell_0}{3}\right)^{\nu}$ (for some $\ell_0\leq\ell$), the distribution of $Y$ is bounded by $\lambda\left(\frac{3}{\ell_0}\right)^{\nu}$ for some positive absolute constant $\nu$. Then an analogous result is implied for the Littlewood-Offord problem over $\mathbb{Z}_p$ on subsets $X$ of a given cardinality $\ell$ in the regime where $n$ is large enough.

Finally, as an application of our results, we propose a variation of the set-sequenceability problem: that of $\Gamma$-sequenceability. Given a graph $\Gamma$ on the vertex set $\{1,2,\dots,n\}$ and given a subset $A\subseteq \mathbb{Z}_p$ of size $n$, here we want to find an ordering of $A$ such that the partial sums $s_i$ and $s_j$ are different whenever $\{i,j\}\in E(\Gamma)$. As a consequence of our results on the Littlewood-Offord problem, we have been able to prove that, if the maximum degree of $\Gamma$ is at most $d$, $n$ is large enough, and $p>n^2$, any subset $A\subseteq \mathbb{Z}_p$ of size $n$ is $\Gamma$-sequenceable.
\end{abstract}
\section{Introduction}
The Littlewood-Offord problem is a classical combinatorial question in probability theory and discrete mathematics that has garnered significant attention in various fields, such as statistics, number theory, and theoretical computer science. I had been proposed, firstly by Littlewood and Offord, in the 1940s (see \cite{LO}) in order to deal with the distribution of sums of random variables with restricted distributions. In its first enunciation, given a list of (not necessarily distinct) integers, the problem seeks to determine the likelihood of obtaining a given number by summing elements from the given list. Equivalently, given a list $(v_1, v_2, \dots, v_{\ell})$ of integers, this problem asks about the probability distribution of the sum $\sum_{i=1}^{\ell} Y_i$ where $Y_i$ is the uniform random variable whose support is $\{0,v_i\}$ (or, in an equivalent formulation, $\{-v_i,v_i\}$). In \cite{LO} Littlewood and Offord proved that this probability is at most $O(\frac{\log{n}}{\sqrt{n}})$, result improved then by Erd\H{o}s to $\frac{1}{\sqrt{n}}(1+o(1))$ in \cite{E}.

More in general, this problem and its variation arise in various settings, such as in the analysis of random walks, the study of random matrices, and the investigation of other combinatorial structures (see, for instance, the paper of Bibak \cite{B} and that of Tao and Vu \cite{TV2}). Significant progress has been done in understanding the behavior of the Littlewood-Offord problem under different conditions, such as considering more general random variables $Y_i$ (see, for instance, the paper of Juskevicius and Kurauskas \cite{JK} and that of Leader and Radcliffe \cite{LR}) or considering the problem over finite groups.
In particular, the problem has been proposed, by Vaughan and Wooley (see \cite{VW}), in the case where the variables $Y_i$ are uniform on the support $\{0,v_i\}$ and interesting bounds have been obtained there, by Griggs in \cite{G}, by Bibak in \cite{B}, and by Juskevicius and Semetulskis in \cite{JS}. In the latter work, the authors have considered also the case where the distribution of the $Y_i$ is upper-bounded by $1/2$.

The investigation of this paper is widely inspired by these works. In particular, we will consider random variables $Y_1,\dots, Y_{\ell}$ whose supports are contained in a cyclic group $\mathbb{Z}_k$ and whose distributions are upper-bounded by $\lambda=1/n$. Note that, if the distributions are all uniform on $A\subseteq \mathbb{Z}_k$, and if $n=|A|$ is large enough, this problem is linked to that of upper-bounding the probability that a subset $X\subseteq A$ of a given cardinality $\ell$ sums to $x$. Also, this problem can be seen as a variation of the classical Littlewood-Offord one where the cardinality of the subset $X$ has been fixed and where the elements $v_1,v_2,\dots,v_n$ are all distinct. We recall that, the case where the $v_i$ are distinct integers has been studied, firstly, by Erd\H{o}s and then by Halaz in \cite{H} (see also the book of Tao and Vu, \cite{TV1}) where a bound of type $O(\frac{1}{n\sqrt{n}})$ has been proved.

The paper is organized as follows.
In Section 2, we will revisit the integer case, presenting a simple proof of a bound (already proved by Juskevicius and Kurauskas in \cite{JK} in a much more general context) and writing explicitly and improving the case $k=3$. Here we consider the sum $Y=Y_1+\dots+Y_{\ell}$ of $\ell$ independent uniformly distributed (over the same set $A$) random variables $Y_i$.
Denoted by $n$ the cardinality of $A$, here we have that there exists an absolute constant $D$ for which $$\left(\max_{x\in \mathbb{Z}} \mathbb{P}[Y=x]\right)\leq\frac{D}{n\sqrt{\ell-1}}.$$
It follows that, however we take $\epsilon$, if $\ell$ is sufficiently large with respect to $\epsilon$ we have that
$$\left(\max_{x\in \mathbb{Z}} \mathbb{P}[Y=x]\right)\leq\frac{\epsilon}{n}.$$ We also prove that, if $\ell=3$, we can take $\epsilon=\frac{3+1/n^2}{4}$ obtaining a non-trivial bound.

Then, in Section 3, we will consider our problem on cyclic groups. Here, if the prime factors of $k$ are large enough (namely larger than $n!$ where $n=|A|$), we can use a Freiman isomorphism of order $\ell$ to obtain an upper-bound in $\mathbb{Z}_k$. To avoid this, not so pleasant, hypothesis on the prime factors, we consider a different approach.
More precisely, in case the order is a prime $p>2 n\left(\frac{\ell_0}{3}\right)^{\nu}$ (for some $\ell_0\leq\ell$), using a discrete Fourier approach, we will be able to prove a similar result to that of the integer case. Here, if $Y=Y_1+\dots+Y_{\ell}$ where $Y_i$ are i.i.d. whose distributions are upper-bounded by $1/n=\lambda\leq 9/10$, $\ell\geq \ell_0$, there exists a positive absolute constant $\nu$ for which $$\left(\max_{x\in \mathbb{Z}_p} \mathbb{P}[Y=x]\right)\leq\frac{1}{n}\left(\frac{3}{\ell_0}\right)^{\nu}.$$

In Section 4 we will come back to the original Littlewood-Offord problem over cyclic groups of prime order $p$. As a consequence of our bounds of Section 3, we can prove that, however, we take $\epsilon$, if $\ell$, $n$, and $n-\ell$ are sufficiently large and if $p>\frac{4n}{\epsilon}$ (or, more simply, $p>n^2$), we have that
$$\left(\max_{x\in \mathbb{Z}_p} \mathbb{P}[\sum_{v\in X} v=x]\right)\leq\frac{\epsilon}{n}$$
where the set $X$ is chosen, uniformly at random, among the subsets of size $\ell$ of $A$ and where the cardinality of $A$ is $n$.

The main application we will propose for these results concerns a variation of the sequenceability problem.
We recall that a subset $A$ of a group $G$ is {\em sequenceable} if there is an ordering $(a_1, \ldots, a_n)$ of its elements such that the partial sums~$(s_0, s_1, \ldots, s_n)$, given by $s_0 = 0$ and $s_i = \sum_{j=1}^i a_j$ for $1 \leq i \leq k$, are distinct, with the possible exception that we may have~$s_n = s_0 = 0$. Several conjectures and questions concerning the sequenceability of subsets of groups arose in the Design Theory context: indeed this problem is related to Heffter Arrays and $G$-regular Graph Decompositions (see the paper of Archdeacon where the Heffter Arrays have been introduced \cite{A15}, the papers \cite{ADMS16, HOS19} about these conjectures, the survey of Ollis \cite{OllisSurvey} and that of Pasotti and Dinitz \cite{PD}). Alspach and Liversidge combined and summarized many of them in \cite{AL20} into the conjecture that, if a subset of an abelian group does not contain $0$, then it is sequenceable.

Note that, if the elements of a sequenceable set $A$ do not sum to $0$, then there exists a simple path $P$ in the Cayley graph $Cay[G:\pm A]$ such that $\Delta(P) = \pm S$. Inspired by this interpretation, Costa and Della Fiore proposed in \cite{CD} a weakening of this concept. In particular, they wanted to find an ordering whose partial sums define a walk $W$, in the Cayley graph $Cay[G:\pm A]$, of girth larger than $t$ (for a given $t < k$) and such that $\Delta(W) = \pm S$. This is possible given that the partial sums $s_i$ and $s_j$ are different whenever $i$ and $j$ are distinct and $|i-j|\leq t$.

The problem we will consider in the last section of this paper is a variation of the previous one. Given a graph $\Gamma$ on the vertex set $\{1,2,\dots,n\}$, here we want to find an ordering of $A$ such that the partial sums $s_i$ and $s_j$ are different whenever $\{i,j\}\in E(\Gamma)$.

In particular, we will be able to prove that, if the maximum degree of $\Gamma$ is bounded by a natural number $d$, $n$ is large enough with respect to $d$ and $p>n^2$, any subset $A$ of cardinality $n$ of $\mathbb{Z}_p\setminus \{0\}$ is $\Gamma$-sequenceable.
\section{The integer case}
In this section, given a set $A=\{v_1,v_2,\dots,v_n\}$ of distinct integers, we consider the random variable $Y$ given by the sum of $\ell$ independent variables $Y_1,\dots, Y_{\ell}$ that are uniformly distributed on $A$. Here we will re-propose, with a direct proof, an upper bound on $\mathbb{P}[Y=x]$ that is a special case of Theorem 2.3 of \cite{JK} (where Juskevicius and Kurauskas considered a much more general setting). First of all we restrict ourself to consider $A=\{-(n-1)/2,\dots,(n-1)/2\}$, indeed Theorem 2 of \cite{LR} state that:
\begin{thm}[Leader and Radcliffe]\label{LRthm} Denoted by $\tilde{Y}$ the random variable given by the sum of $\ell$ independent and uniformly distributed on $\{-(n-1)/2,\dots,(n-1)/2\}$ variables $\tilde{Y}_1,\dots, \tilde{Y}_{\ell}$ and considered $Y$ defined as above, we have that
$$\left(\max_{x\in \mathbb{Z}} \mathbb{P}[Y=x]\right)\leq\left(\max_{x\in \mathbb{Z}} \mathbb{P}[\tilde{Y}=x]\right)=\mathbb{P}[\tilde{Y}=M]$$
where $M=\begin{cases}0 \mbox{ if } \ell\equiv 0 \pmod{2} \mbox{ or }n\equiv 1 \pmod{2};\\
-\frac{1}{2} \mbox{ otherwise}.\end{cases}$
\end{thm}
Then our main tool to upper-bound $\mathbb{P}[\tilde{Y}=x]$ is the Berry-Esseen theorem which says that:
\begin{thm}[Berry-Esseen]\label{BE}
Let $Y_1,\dots, Y_{\ell'}$ be independent random variables with the same distribution and such that:
\begin{itemize}
\item[(a)] $\mathbb{E}(Y_j)=0$;
\item[(b)] $\mathbb{E}(Y_j^2)=\sigma^2>0$;
\item[(c)] $\mathbb{E}(|Y_j|^3)<\infty$.
\end{itemize}
Then, set $\rho=\frac{\mathbb{E}(|Y_j|^3)}{\sigma^3}$ we have

$$\sup_x \left| \mathbb{P}[\frac{1}{\sigma\sqrt{\ell'}}\sum_{j=1}^{\ell'} Y_j\leq x]-\Phi(x) \right|\leq C\frac{\rho}{\sqrt{\ell'}} $$
where $\Phi(x)$ is the cumulative distribution function of the standard normal distribution and $C\leq 0.7655$ (see \cite{S}) is an absolute constant.
\end{thm}
Then we can state the following:
\begin{thm}\label{Z,ell,epsilon}
Let $A=\{v_1,v_2,\dots,v_n\}$ be distinct integers and let $Y$ be the sum of $\ell$ independent and uniformly distributed on $A$ variables $Y_1,\dots, Y_{\ell}$.
Then, there exists an absolute constant $D$ for which $$\left(\max_{x\in \mathbb{Z}} \mathbb{P}[Y=x]\right)\leq\frac{D}{n\sqrt{\ell-1}}.$$
In particular, however we take $\epsilon$, if $\ell$ is sufficiently large with respect to $\epsilon$ we have that
$$\left(\max_{x\in \mathbb{Z}} \mathbb{P}[Y=x]\right)\leq\frac{\epsilon}{n}.$$
\end{thm}
\proof
Due to Theorem \ref{LRthm}, we can assume $A=\{-(n-1)/2,\dots,(n-1)/2\}$.
We note that, due to the independence of the $Y_i$,
$$\mathbb{P}[Y=x]=$$
\begin{equation}\label{riduzione}\sum_{y\in [x-\frac{n-1}{2},x+\frac{n-1}{2}]}\mathbb{P}[Y_1+\dots+Y_{\ell-1}=y]\mathbb{P}[Y_{\ell}=x-y]=\frac{\mathbb{P}[Y_1+\dots+Y_{\ell-1}\in [x-\frac{n-1}{2},x+\frac{n-1}{2}]]}{n}.\end{equation}
Now we want to extimate $\mathbb{P}[Y_1+\dots+Y_{\ell-1}\in [x-\frac{n-1}{2},x+\frac{n-1}{2}]]$ through Theorem \ref{BE}.
Indeed, set
$$\Psi(z)=\mathbb{P}[\frac{1}{\sigma\sqrt{\ell-1}}\sum_{j=1}^{\ell-1} Y_j\leq z]$$
and set $(\Psi-\Phi)(z):=\Psi(z)-\Phi(z)$, we have:
$$\mathbb{P}[Y_1+\dots+Y_{\ell-1}\in [x-\frac{n-1}{2},x+\frac{n-1}{2}]]=\Psi\left(\frac{(x+\frac{n-1}{2})}{\sigma\sqrt{\ell-1}}\right)-\Psi\left(\frac{(x-\frac{n-3}{2})}{\sigma\sqrt{\ell-1}}\right)$$
\begin{equation}\label{Spezzata1}=\left(\Psi-\Phi\right)\left(\frac{(x+\frac{n-1}{2})}{\sigma\sqrt{\ell-1}}\right)-\left(\Psi-\Phi\right)\left(\frac{(x-\frac{n-3}{2})}{\sigma\sqrt{\ell-1}}\right)+\end{equation}
$$\Phi\left(\frac{(x+\frac{n-1}{2})}{\sigma\sqrt{\ell-1}}\right)-\Phi\left(\frac{(x-\frac{n-3}{2})}{\sigma\sqrt{\ell-1}}\right).$$
This implies that, due to the triangular inequality,
$$\mathbb{P}[Y_1+\dots+Y_{\ell-1}\in [x-\frac{n-1}{2},x+\frac{n-1}{2}]]\leq $$
\begin{equation}\label{Spezzata2}\left|\left(\Psi-\Phi\right)\left(\frac{(x+\frac{n-1}{2})}{\sigma\sqrt{\ell-1}}\right)\right|+\left|\left(\Psi-\Phi\right)\left(\frac{(x-\frac{n-3}{2})}{\sigma\sqrt{\ell-1}}\right)\right|+\end{equation}
$$\left|\Phi\left(\frac{(x+\frac{n-1}{2})}{\sigma\sqrt{\ell-1}}\right)-\Phi\left(\frac{(x-\frac{n-3}{2})}{\sigma\sqrt{\ell-1}}\right)\right|.$$
We note that the first two terms of Equation \eqref{Spezzata2} can be upper bounded via the Berry-Esseen theorem.
Here we have that
$$\sigma=\frac{n}{2\sqrt{3}}(1+o(1))$$
and
$$\rho=\frac{3\sqrt{3}}{4}(1+o(1)).$$
Therefore
\begin{equation}\label{primidue}\left|\left(\Psi-\Phi\right)\left(\frac{(x+\frac{n-1}{2})}{\sigma\sqrt{\ell-1}}\right)\right|+\left|\left(\Psi-\Phi\right)\left(\frac{(x-\frac{n-3}{2})}{\sigma\sqrt{\ell-1}}\right)\right|\leq 2C\frac{3\sqrt{3}}{4\sqrt{\ell-1}}(1+o(1)).\end{equation}
For the latter term, we note that
$$\left|\Phi\left(\frac{(x+\frac{n-1}{2})}{\sigma\sqrt{\ell-1}}\right)-\Phi\left(\frac{(x-\frac{n-3}{2})}{\sigma\sqrt{\ell-1}}\right)\right|\leq\max_x |\Phi(x)| \frac{n}{\sigma\sqrt{\ell-1}}.$$
Hence, since the cumulative distribution function of the standard normal distribution has maximum $\frac{1}{2\pi}$ in $0$, we have
\begin{equation}\label{terzotermine}
\left|\Phi\left(\frac{(x+\frac{n-1}{2})}{\sigma\sqrt{\ell-1}}\right)-\Phi\left(\frac{(x-\frac{n-3}{2})}{\sigma\sqrt{\ell-1}}\right)\right|\leq \frac{1}{2\pi}\frac{2\sqrt{3}}{\sqrt{l-1}}(1+o(1)).
\end{equation}
Summing up Equations \eqref{primidue} and \eqref{terzotermine}, we have that
$$\mathbb{P}[Y_1+\dots+Y_{\ell-1}\in [x-\frac{n-1}{2},x+\frac{n-1}{2}]]\leq \frac{D}{\sqrt{\ell-1}}$$
where $D$ is an absolute constant that also incorporates the term $o(1)$ which is bounded and goes to zero when $n$ goes to infinite.
Placing this bound in Equation \eqref{riduzione}, we obtain that
$$\mathbb{P}[Y=x]\leq \frac{D}{n\sqrt{\ell-1}}$$
which implies the thesis.
\endproof
Even though for $\ell=3$ the previous theorem provides a trivial bound, we can deal with this case directly. More precisely we obtain the following non-trivial bound:
\begin{thm}\label{Z,3,epsilon}
Let $A=\{v_1,v_2,\dots,v_n\}$ be distinct integers, and let $Y_1, Y_2, Y_3$ be independent and uniformly distributed on $A$. Then, set $Y=Y_1+Y_2+Y_3$, we have that
$$\left(\max_{x\in \mathbb{Z}} \mathbb{P}[Y=x]\right)\leq\frac{3+1/n^2}{4n}.$$
\end{thm}
\proof
Due to Theorem \ref{LRthm}, we can assume $A=\{1,\dots,n\}$ (which is a traslate of $\{-(n-1)/2,\dots,(n-1)/2\}$).
Also, Theorem \ref{LRthm} says that
$\mathbb{P}[Y=x]$ is maximal for $x=M=\lfloor\frac{3(n+1)}{2}\rfloor$ when $A=\{1,\dots,n\}$ (which corresponds to $x=0$ or $-\frac{1}{2}$ for $\{-(n-1)/2,\dots,(n-1)/2\}$).
Now we compute the number of triples
$$T=\{(y_1,y_2,y_3)\in A^3: y_1+y_2+y_3=M\}.$$

We consider explicitly only the case $n$ odd: the case $n$ even is completely analogous. Clearly $y_1+y_2+y_3=M$ implies that $y_1+y_2\in [M-n,M-1]=[\frac{n+3}{2},\frac{3n+1}{2}]$ and that $y_3=M-y_2-y_1$.

Given $x\in [\frac{n+3}{2},\frac{3n+1}{2}]$, we need to find the number of pair $(y_1,y_2)$ which sum to $x$. Here we have that, for any $y_2$ such that $1\leq x-y_2\leq n$, there exists exactly one $y_1$ such that $y_1+y_2=x$.
It follows that, if $x\leq n$ there are exactly $x-1$ pairs $y_1,y_2$ which sums to $x$. Similarly, if $x\geq n+1$, there are exactly $2n-(x-1)$ pairs $y_1,y_2$ which sums to $x$.
Noting that the number of pairs that sum to $x$ is equal to that of the pairs which sum to $2n+2-x$ we obtain that:
$$
|T|=\sum_{x\in [\frac{n+3}{2},n]}(x-1)+\sum_{x\in [n+1,\frac{3n+1}{2}]}(2n-(x-1))=
$$
\begin{equation}\label{numbertriple}2\sum_{x\in [\frac{n+3}{2},n]}(x-1)+n=2\left(\sum_{x=1}^{n-1} x-\sum_{x=1}^{\frac{n-1}{2}}x\right)+n=\frac{3n^2+1}{4}.\end{equation}
The thesis follows because
$$\mathbb{P}[Y=x]= \frac{\{(y_1,y_2,y_3)\in A^3: y_1+y_2+y_3=x\}}{n^3}\leq \frac{T}{n^3}\leq \frac{3+1/n^2}{4n}.$$
\endproof
\section{Anticoncentration inequalities for $\mathbb{Z}_k$ and $\mathbb{Z}_p$}
In Section 2 we have presented an upper bound on $\mathbb{P}[Y=x]$ for $\mathbb{Z}$. It is then natural to investigate the case of $\mathbb{Z}_k$. Here, given a subset $A$ of size $n$, if the prime factors of $k$ are large enough, there exists a Freiman isomorphism of order $\ell$ from $A$ to a subset $B\subseteq \mathbb{Z}$. We recall that (see Tao and Vu, \cite{TV1}).
\begin{defn}
Let $k\geq 1$, and let $A$, $B$ be additive sets with ambient groups $V$ and $W$ respectively. A Freiman homomorphism of order
$\ell$, say $\phi$, from $(A, V)$ to $(B,W)$ (or more succinctly from $A$ to $B$) is a map $\phi : A \rightarrow B$ with the property that:
$$a_1+a_2+\dots+a_{\ell}=a_1'+a_2'+\dots+a_{\ell}'\rightarrow \phi(a_1)+\phi(a_2)+\dots+\phi(a_{\ell})=\phi(a_1')+\phi(a_2')+\dots+\phi(a_{\ell}').$$
If in addition there is an inverse map $\phi^{-1}: B \rightarrow A$ which is also a Freiman homomorphism of order $\ell$ from $(B,W)$ to $(A, V)$, then we say that $\phi$ is a Freiman isomorphism of order $\ell$ and that $(A, V)$ and $(B,W)$ are Freiman isomorphic of order $\ell$.
\end{defn}
Also, following the proof of Theorem 22 of \cite{CDOR} (see also Theorem 3.1 of \cite{CP20} where an explicit bound was not provided), we have that
\begin{thm}[\cite{CDOR}]
Let $A$ be a subset of size $n$ of $\mathbb{Z}_k$ where the prime factors of $k$ are larger than $n!$. Then there exists a Freiman isomorphism of any order $\ell\leq n$, from $A$ to a subset $B\subseteq \mathbb{Z}$.
\end{thm}
Therefore, the result over the integers, namely Theorem \ref{Z,ell,epsilon}, implies that
\begin{thm}\label{Zk,ell,epsilon}
Let $A$ be a subset of size $n$ of $\mathbb{Z}_k$ where the prime factors of $k$ are larger than $n!$. Then, set $Y=Y_1+\dots+Y_{\ell}$ the sum of $\ell$ independent uniformly distributed (over the same set $A$) random variables $Y_i$,
$$\mathbb{P}[Y=x]\leq \frac{D}{n\sqrt{\ell-1}}$$
where $D$ is an absolute constant.
\end{thm}
Note that this bound has a quite strong hypothesis on the prime factors of $k$. Therefore, in following of this section, we consider a different approach to improve it. Now we consider a set $A=\{v_1,v_2,\dots,v_n\}$ of distinct elements of $\mathbb{Z}_p$, and the random variable $Y$ given by the sum of $\ell$ independent and uniformly distributed on $A$ variables $Y_1,\dots, Y_{\ell}$. The goal is here to provide an upperbound on $\left(\max_{x\in \mathbb{Z}_p} \mathbb{P}[Y=x]\right)$. First of all, we consider the case $\ell=3$ and we prove the following Lemma.
\begin{lem}\label{UniformSym}
Let $n\geq 2$, $p>2n$, $A=\{v_1,v_2,\dots,v_n\}$ be a set of distinct elements of $\mathbb{Z}_p$ such that $-v\in A$ whenever $v\in A$ and let $Y_1,Y_2$ and $Y_3$ be independent variables which are uniformly distributed on $A$. Then there exists an absolute constant $C_1<1$\footnote{Here the best approximation we have for this constant is $C_1< 0.99993.$}, such that, set $Y=Y_1+Y_2+Y_3,$ we have:
$$\left(\max_{x\in \mathbb{Z}_p} \mathbb{P}[Y=x]\right)\leq\frac{C_1}{n}.$$
\end{lem}
\proof
First of all, we note that, if $n=2$, we may suppose, without loss of generality, that $A=\{-1,1\}$. For this set, with the same proof of Theorem \ref{Z,3,epsilon}, we obtain that $\mathbb{P}[Y=x]\leq \frac{3+1/4}{4n}$ (which give a better constant than $\frac{C_1}{n}$). So, in the following, we only consider the case $n\geq 3$.

Following the proof of Proposition 6.1 of \cite{FKS}, set $f=\sum_{i=1}^n \frac{1}{n}\delta_{v_i}$, we have that its discrete Fourier transform is
$$ \hat{f}(k)=\sum_{i=1}^n \frac{1}{n}e^{-2\pi i v_ik/p}.$$
Then, the probability that $Y_1+Y_2=x$ can be written as
$$\mathbb{P}[Y_1+Y_2=x]=(f*f)(x)$$
$$=\frac{1}{p}\sum_{k=0}^{p-1} e^{2\pi ixk/p} (\hat{f}(k))^2$$
$$=\frac{1}{p}\sum_{k=0}^{p-1} e^{2\pi ixk/p} \left(\sum_{i=1}^n \frac{1}{n}e^{-2\pi i v_ik/p}\right)^2.$$
Here we note that, since the set $A$ is symmetric with respect to $0$, $$\mathbb{P}[Y_1+Y_2=x]=\mathbb{P}[Y_1+Y_2=-x].$$
Thus we can write
\begin{equation}\label{cos1}
\mathbb{P}[Y_1+Y_2=x]=\frac{1}{p}\sum_{k=0}^{p-1} \cos( 2\pi xk/p) \left(\sum_{i=1}^n \frac{1}{n}e^{-2\pi i v_ik/p}\right)^2.
\end{equation}
Since, due again to the symmetry of $A$, we also have that both $e^{-2\pi i v_ik/p}$ and $e^{2\pi i v_ik/p}$ appears in $\hat{f}(k)$, Equation \eqref{cos1} can be written as:
\begin{equation}\label{cos2}
\mathbb{P}[Y_1+Y_2=x]=\frac{1}{p}\sum_{k=0}^{p-1} \cos( 2\pi xk/p) \left(\sum_{i=1}^n \frac{1}{n}\cos(2\pi v_ik/p)\right)^2.
\end{equation}
Now we consider $B\subseteq \mathbb{Z}_p$ of cardinality $n$ such that
$\mathbb{P}[Y_1+Y_2\in B]$ is maximal. We divide $B$ into two sets defined as follows:
$$B_1:=\{x\in B: \mathbb{P}[Y_1+Y_2=x]\geq \frac{1-\epsilon_1}{n}\};$$
$$B_2:=B\setminus B_1=\{x\in B: \mathbb{P}[Y_1+Y_2=x]< \frac{1-\epsilon_1}{n}\}.$$
Now we divide the proof into two cases.

CASE 1: $|B_1|\geq (1-\epsilon_2)n$. Here we note that
\begin{equation}\label{val0}
\mathbb{P}[Y_1+Y_2=0]=\frac{1}{n}=\frac{1}{p}\sum_{k=0}^{p-1}(\hat{f}(k))^2.
\end{equation}
This means that, for the values $x \in B_1$, we can not lose too much by inserting the cosine. Indeed, for these values, we have that
\begin{equation}\label{valx}
\mathbb{P}[Y_1+Y_2=x]=\frac{1}{p}\sum_{k=0}^{p-1} \cos( 2\pi xk/p)(\hat{f}(k))^2\geq \frac{1-\epsilon_1}{n}.
\end{equation}
Taking the average over the values $x\in B_1$, we obtain
$$\frac{1}{|B_1|}\mathbb{P}[Y_1+Y_2\in B_1]=\frac{1}{|B_1|}\sum_{x\in B_1}\frac{1}{p}\sum_{k=0}^{p-1} \cos( 2\pi xk/p)(\hat{f}(k))^2=
$$
\begin{equation}\label{average}
=\frac{1}{p}\sum_{k=0}^{p-1} \left ((\hat{f}(k))^2 \frac{\sum_{x\in B_1}\cos( 2\pi xk/p)}{|B_1|}\right)=\frac{1}{p}+\frac{1}{p}\sum_{k=1}^{p-1} \left ((\hat{f}(k))^2 \frac{\sum_{x\in B_1}\cos( 2\pi xk/p)}{|B_1|}\right),
\end{equation}
where the last equalities hold by commuting the sums and noting that $\hat{f}(0)=1$.
Now, Equations \eqref{val0} and \eqref{valx} implies that
$$\frac{1}{|B_1|}\mathbb{P}[Y_1+Y_2\in B_1] \geq \frac{1-\epsilon_1}{n}=\frac{1}{p}\sum_{k=0}^{p-1}\left((\hat{f}(k))^2(1-\epsilon_1)\right)=$$ $$=\frac{1-\epsilon_1}{p}+\frac{1}{p}\sum_{k=1}^{p-1}\left((\hat{f}(k))^2(1-\epsilon_1)\right).$$
Now since $2n< p$, Equation \eqref{average} implies the existence of $\bar{k}\not=0$ such that \begin{equation}\label{bark}\sum_{x\in B_1}\frac{\cos( 2\pi x\bar{k}/p)}{|B_1|}\geq (1-2\epsilon_1).\end{equation}
Indeed, otherwise, we would have that
$$\frac{1}{p}+\frac{1}{p}\sum_{k=1}^{p-1} \left ((\hat{f}(k))^2 (1-2\epsilon_1\right)> \frac{1}{p}+\frac{1}{p}\sum_{k=1}^{p-1} \left ((\hat{f}(k))^2 \frac{\sum_{x\in B_1}\cos( 2\pi xk/p)}{|B_1|}\right)\geq $$ $$ \frac{1-\epsilon_1}{p}+\frac{1}{p}\sum_{k=1}^{p-1}\left((\hat{f}(k))^2(1-\epsilon_1)\right)$$
that is, since $2n\leq p$,
$$\frac{\epsilon_1}{p}> \frac{1}{p}\sum_{k=1}^{p-1}(\hat{f}(k))^2\epsilon_1=\epsilon_1\left(\frac{1}{n}-\frac{1}{p} \right)\geq \frac{\epsilon_1}{p}$$
which is a contradiction. Also, note that since $\mathbb{Z}_p$ is a field, we may suppose, without loss of generality, that $\bar{k}=1$.

This means that
\begin{equation}\label{cos3}\sum_{x\in B_1}\frac{\cos( 2\pi x/p)}{|B_1|}\geq (1-2\epsilon_1).\end{equation}
Hence the set
$$B_1':=\{x\in B_1: x/p\in [-1/6,1/6]\}=\{x\in B_1: \cos( 2\pi x/p)\in [1/2,1]\}$$ has size larger than $(1-2\epsilon_2-4\epsilon_1)n$.
Indeed, since $|B|\leq n$, $|B_1'|<(1-2\epsilon_2-4\epsilon_1)n$ would imply that
$$\frac{n+(1-2\epsilon_2-4\epsilon_1)n}{2} \geq \frac{|B_1|+|B_1'|}{2}=|B_1'|+\frac{|B_1|-|B_1'|}{2}.$$
Then we would obtain
$$(1-\epsilon_2-2\epsilon_1)n>|B_1'|+\frac{|B_1\setminus B_1'|}{2}\geq \sum_{x\in B_1}\cos( 2\pi x/p)\geq (1-2\epsilon_1)(1-\epsilon_2)n$$
which is a contradiction.

Now we come back to the estimation of $\max_{x\in \mathbb{Z}_p}\mathbb{P}[Y=x]$. Because of the independence of the variables $Y_1,Y_2$ and $Y_3$, we have
\begin{equation}\label{independence}
\mathbb{P}[Y=x]=\sum_{v_i\in A}\mathbb{P}[Y_1+Y_2=x-v_i]\mathbb{P}[Y_3=v_i].
\end{equation}
Now we split this computation according to whether $x-v_i$ belongs to $B_1'$ to $B_1'':=B_1\setminus B_1'$ or $B_2$.
\begin{equation}\label{split1}
\sum_{v_i\in A:\ x-v_i\in B_1'}\mathbb{P}[Y_1+Y_2=x-v_i]\mathbb{P}[Y_3=v_i]+\sum_{v_i\in A:\ x-v_i\in B_1''}\mathbb{P}[Y_1+Y_2=x-v_i]\mathbb{P}[Y_3=v_i]+
\end{equation}
$$+\sum_{v_i\in A:\ x-v_i\in B_2}\mathbb{P}[Y_1+Y_2=x-v_i]\mathbb{P}[Y_3=v_i].$$
Due to Equations \eqref{val0} and \eqref{valx}, $\mathbb{P}[Y_1+Y_2=x-v_i]$ is always bounded by $1/n$. Then, since $|\{v_i\in A:\ x-v_i\in B_1''\}|\leq (4\epsilon_1+2\epsilon_2)n$, we can provide the following bounds:
$$\sum_{v_i\in A:\ x-v_i\in B_1'}\mathbb{P}[Y_1+Y_2=x-v_i]\mathbb{P}[Y_3=v_i]\leq\frac{1}{n} \frac{|\{v_i\in A:\ x-v_i\in B_1'\}|}{n};$$
$$\sum_{v_i\in A:\ x-v_i\in B_1''}\mathbb{P}[Y_1+Y_2=x-v_i]\mathbb{P}[Y_3=v_i]\leq \frac{4\epsilon_1+2\epsilon_2}{n};$$
and
$$\sum_{v_i\in A:\ x-v_i\in B_2}\mathbb{P}[Y_1+Y_2=x-v_i]\mathbb{P}[Y_3=v_i]<\frac{1-\epsilon_1}{n}\frac{|\{v_i\in A:\ x-v_i\in B_2\}|}{n}\leq \frac{\epsilon_2(1-\epsilon_1)}{n}.$$
Therefore $\mathbb{P}[Y=\bar{x}]\geq\frac{1-\epsilon_2}{n}+\frac{\epsilon_2(1-\epsilon_1)}{n}$ (the right-hand side of this inequality will be denoted in the following by $\frac{C_1}{n}$) implies that
$$\frac{1}{n} \frac{|\{v_i\in A:\ \bar{x}-v_i\in B_1'\}|}{n}+\frac{4\epsilon_1+2\epsilon_2}{n}+\frac{\epsilon_2(1-\epsilon_1)}{n^2}\geq \frac{1-\epsilon_2}{n}+\frac{\epsilon_2(1-\epsilon_1)}{n}.$$
As a consequence, we have
$$|\{v_i\in A:\ \bar{x}-v_i\in B_1'\}|\geq (1-4\epsilon_1-3\epsilon_2)n$$
and hence, set
$$A_1:=\{v_i\in A:\ v_i\in [-1/6-\bar{x},1/6-\bar{x}]\},$$
we have that
\begin{equation}\label{vi}
|A_1|\geq |\{v_i\in A:\ \bar{x}-v_i\in B_1'\}|\geq (1-4\epsilon_1-3\epsilon_2)n.
\end{equation}
On the other hand, given two triples $v_1,v_2,v_3$ and $v_1',v_2',v_3'\in A_1$, we have that $v_1+v_2+v_3=v_1'+v_2'+v_3'$ if and only if have the same sum also in $\mathbb{Z}$ (with the trivial identification).
Hence, due to Theorem \ref{Z,3,epsilon} and since $n\geq 3$,
$$\max_{x\in \mathbb{Z}_p} \mathbb{P}[Y=x|Y_1,Y_2,Y_3\in A_1]\leq \frac{3+1/9}{4n (1-4\epsilon_1-3\epsilon_2)}.$$
Here we have that
$$\mathbb{P}[Y=x|Y_1\not\in A_1]=\mathbb{P}[Y=x|Y_2\not\in A_1]=\mathbb{P}[Y=x|Y_3\not\in A_1],$$ 
and that $\mathbb{P}[Y=x|Y_1\not\in A_1]$ can be written as
$$\sum_{y\in \mathbb{Z}_p}\mathbb{P}[Y_1+Y_2=y|Y_1\not\in A_1]\mathbb{P}[Y_3=x-y]\leq \sum_{y\in \mathbb{Z}_p}\mathbb{P}[Y_1+Y_2=y|Y_1\not\in A_1] \frac{1}{n}=\frac{1}{n}.$$
Hence we can bound $\mathbb{P}[Y=x]$ as follows:
$$\mathbb{P}[Y=x]\leq \mathbb{P}[Y=x|Y_1,Y_2,Y_3\in A_1]+3\mathbb{P}[Y=x|Y_1\not\in A_1]\mathbb{P}[Y_1\not\in A_1]\leq$$
\begin{equation}\label{case1}
\frac{3+1/9}{4n (1-4\epsilon_1-3\epsilon_2)}+\frac{3}{n}\frac{|A\setminus A_1|}{n}\leq \frac{3+1/9}{4n (1-4\epsilon_1-3\epsilon_2)} +3\frac{4\epsilon_1+3\epsilon_2}{n}.
\end{equation}
If we consider $(\epsilon_1,\epsilon_2)$ on of the surface
$$\frac{3+1/9}{4n (1-4\epsilon_1-3\epsilon_2)} +3\frac{4\epsilon_1+3\epsilon_2}{n}=\frac{1-\epsilon_2}{n}+\frac{\epsilon_2(1-\epsilon_1)}{n},$$
we obtain that
$$\mathbb{P}[Y=x]\leq \frac{1-\epsilon_2}{n}+\frac{\epsilon_2(1-\epsilon_1)}{n}=\frac{C_1}{n}$$
and we conclude CASE 1 by computing (with Mathematica) the minimum possible $C_1$ in this surface which is $C_1< 0.99993$.

CASE 2: $|B_1|< (1-\epsilon_2)n$. Here we note that Equation \eqref{independence} can be split as follows:
$$\mathbb{P}[Y=x]=$$
\begin{equation}\label{split2}
\sum_{v_i\in A:\ x-v_i\in B_1}\mathbb{P}[Y_1+Y_2=x-v_i]\mathbb{P}[Y_3=v_i]+\sum_{v_i\in A:\ x-v_i\in B_2}\mathbb{P}[Y_1+Y_2=x-v_i]\mathbb{P}[Y_3=v_i].
\end{equation}
Now we can provide the following bounds:
$$\sum_{v_i\in A:\ x-v_i\in B_1}\mathbb{P}[Y_1+Y_2=x-v_i]\mathbb{P}[Y_3=v_i]\leq \frac{|\{v_i\in A:\ x-v_i\in B_1\}|}{n}\frac{1}{n}$$
and
$$\sum_{v_i\in A:\ x-v_i\in B_2}\mathbb{P}[Y_1+Y_2=x-v_i]\mathbb{P}[Y_3=v_i]\leq \frac{|\{v_i\in A:\ x-v_i\in B_2\}|(1-\epsilon_1)}{n}\frac{1}{n}.$$
Therefore, since $|B_1|< (1-\epsilon_2)n$, we upper-bound the right hand side of Equation \eqref{split2} by assuming that $|\{v_i\in A:\ x-v_i\in B_1\}|= (1-\epsilon_2)n$ and
$|\{v_i\in A:\ x-v_i\in B_2\}|=\epsilon_2 n.$
Summing up, we obtain that
\begin{equation}\label{case2}
\mathbb{P}[Y=x]\leq \frac{1-\epsilon_2}{n}+\frac{\epsilon_2(1-\epsilon_1)}{n}=\frac{C_1}{n}
\end{equation}
which concludes CASE 2.
\endproof

\begin{prop}\label{Uniform}
Let $n\geq 2$, $p>2n$, $A=\{v_1,v_2,\dots,v_n\}$ be distinct elements of $\mathbb{Z}_p$ and let $Y_1,Y_2$ and $Y_3$ be independent variables which are uniformly distributed on $A$. Then there exists $C_2<1$\footnote{Here the best approximation we have for this constant is $C_2< 0.999986.$}, such that, set $Y=Y_1+Y_2+Y_3$, we have:
$$\left(\max_{x\in \mathbb{Z}_p} \mathbb{P}[Y=x]\right)\leq\frac{C_2}{n}.$$
\end{prop}
\proof
First of all, we note that, if $n=2$, we may suppose, without loss of generality, that $A=\{-1,1\}$. For this set, with the same proof of Theorem \ref{Z,3,epsilon}, we obtain that $\mathbb{P}[Y=x]\leq \frac{3+1/4}{4n}$ (which give a better constant than $\frac{C_2}{n}$). So, in the following, we only consider the case $n\geq 3$.

Let $\bar{x}$ be such that $\mathbb{P}[Y_1+Y_2=\bar{x}]$ is maximal. Since $\mathbb{Z}_p$ is a field, we may suppose, without loss of generality, that $\bar{x}=0$.
Set $\epsilon_3$ such that
$$\frac{C_1}{(1-\epsilon_{3})}+3\epsilon_{3}=1-\epsilon_{3}$$ and $C_2=1-\epsilon_3$, we note that $C_2<0.999986$. Then we divide the proof into two cases.

CASE 1: $\mathbb{P}[Y_1+Y_2=0]\geq \frac{1-\epsilon_{3}}{n}$.
Here we define
$$A_1:=\{x\in A:\ -x\in A\}.$$
We have that
$$\mathbb{P}[Y_1+Y_2=0]=\mathbb{P}[Y_1\in A_1]\cdot \mathbb{P}[Y_2=-Y_1]=\frac{|A_1|}{n^2}.$$
This implies that
\begin{equation}\label{sizeA1}|A_1|\geq n(1-\epsilon_{3}).\end{equation}
Note that, since $n\geq 3$, $|A_1|\geq 2$. Hence, due to Lemma \ref{UniformSym}, there exists $C_1$ such that
$$\mathbb{P}[Y=x|Y_1,Y_2,Y_3\in A_1]\leq \frac{C_1}{n(1-\epsilon_{3})}.$$
Here we have that
$$\mathbb{P}[Y=x|Y_1\not\in A_1]=\mathbb{P}[Y=x|Y_2\not\in A_1]=\mathbb{P}[Y=x|Y_3\not\in A_1],$$
and that $\mathbb{P}[Y=x|Y_1\not\in A_1]$ can be written as
$$\sum_{y\in \mathbb{Z}_p}\mathbb{P}[Y_1+Y_2=y|Y_1\not\in A_1]\mathbb{P}[Y_3=x-y]\leq \sum_{y\in \mathbb{Z}_p}\mathbb{P}[Y_1+Y_2=y|Y_1\not\in A_1] \frac{1}{n}=\frac{1}{n}.$$
Therefore, we can bound $\mathbb{P}[Y=x]$ as follows:
$$\mathbb{P}[Y=x]\leq \mathbb{P}[Y=x|Y_1,Y_2,Y_3\in A_1]+3\mathbb{P}[Y=x|Y_1\not\in A_1]\mathbb{P}[Y_1\not\in A_1]\leq$$
\begin{equation}\label{case1b}
\frac{C_1}{n(1-\epsilon_{3})}+\frac{3}{n}\frac{|A\setminus A_1|}{n}\leq \frac{C_1}{n(1-\epsilon_{3})}+3\frac{\epsilon_{3}}{n}=\frac{C_2}{n}
\end{equation}
which concludes CASE 1.

CASE 2: $\mathbb{P}[Y_1+Y_2=0]< \frac{1-\epsilon_{3}}{n}$. Here we have that
\begin{equation}\label{split3}
\mathbb{P}[Y=x]=\sum_{v_i\in A}\mathbb{P}[Y_1+Y_2=x-v_i]\mathbb{P}[Y_3=v_i].
\end{equation}
Since
$$\mathbb{P}[Y_1+Y_2=x-v_i]\leq \mathbb{P}[Y_1+Y_2=0]<\frac{1-\epsilon_{3}}{n},$$ Equation \eqref{split3} can be written as
\begin{equation}\label{case2b}
\mathbb{P}[Y=x]<\sum_{v_i\in A}\frac{1-\epsilon_{3}}{n}\mathbb{P}[Y_3=v_i]=\frac{1-\epsilon_{3}}{n}=\frac{C_2}{n}.
\end{equation}
\endproof
This proposition can be generalized to:
\begin{thm}\label{General3}
Let $\lambda\leq \frac{9}{10}$, $p> \frac{2}{\lambda}$, $Y_1,Y_2$ and $Y_3$ be identical and independent variables such that $$\left(\max_{x\in \mathbb{Z}_p} \mathbb{P}[Y_1=x]\right)\leq \lambda.$$ Then there exists $C_3<1$\footnote{Here the best approximation we have for this constant is $C_3< 1-2.27\cdot 10^{-12}.$}, such that, set $Y=Y_1+Y_2+Y_3$ we have:
$$\left(\max_{x\in \mathbb{Z}_p} \mathbb{P}[Y=x]\right)\leq C_3\lambda.$$
\end{thm}
\proof
We consider $\epsilon_{4}$ and $\epsilon_5$ in the surface:
$$1-\epsilon_4\epsilon_5=3(\epsilon_4+\epsilon_5)+\frac{C_2}{(1-\epsilon_4)^3(1-\epsilon_5)}$$
and we set $C_3:=1-\epsilon_4\epsilon_5$. We will also assume that $\epsilon_5<1/10$.

Then, we define
$$A_1:=\{x\in \mathbb{Z}_p: \mathbb{P}[Y_1=x]\geq \lambda(1-\epsilon_{4})\}.$$
Then we divide the proof into two cases according to the cardinality of $A_1$.

CASE 1: $|A_1|\geq \frac{(1-\epsilon_{5})}{\lambda}.$
We want to provide an upper bound to $\mathbb{P}[Y=x|Y_1,Y_2,Y_3\in A_1]$ by applying Proposition \ref{Uniform}. For this purpose, we note that:
$$\mathbb{P}[Y=x|Y_1,Y_2,Y_3\in A_1]=$$
\begin{equation}\label{Stima1}
\sum_{y_1,y_2\in A_1}\mathbb{P}[Y_1=y_1|Y_1\in A_1]\mathbb{P}[Y_2=y_2|Y_2\in A_1]\mathbb{P}[Y_3=x-y_1-y_2|Y_3\in A_1].
\end{equation}
Here, for $y_1\in A_1$, we have
$$
\mathbb{P}[Y_1=y_1|Y_1\in A_1]\mathbb{P}[Y_1\in A_1]=\mathbb{P}[Y_1=y_1]\leq \lambda.
$$
Since
$$\mathbb{P}[Y_1\in A_1]=\sum_{x\in A_1} \mathbb{P}[Y_1=x]\geq |A_1| \lambda (1-\epsilon_{4}),$$
it follows that
$$
\mathbb{P}[Y_1=y_1|Y_1\in A_1]\leq \frac{\lambda}{\mathbb{P}[Y_1\in A_1]}\leq \frac{1}{|A_1|(1-\epsilon_{4})}.
$$
Note that, named by $\tilde{Y}=\tilde{Y_1}+\tilde{Y_2}+\tilde{Y_3}$ where $\tilde{Y_1}, \tilde{Y_2}, \tilde{Y_3}$ are uniform distribution over $A_1$, we have that
\begin{equation}\label{Bound3}
\mathbb{P}[Y_1=y_1|Y_1\in A_1] \leq \frac{1}{|A_1|(1-\epsilon_4)}=\frac{1}{(1-\epsilon_4)}\mathbb{P}[\tilde{Y_1}=y_1].\end{equation}
Noting that $\mathbb{P}[\tilde{Y_1}=y_1]\leq \frac{\lambda}{(1-\epsilon_5)}$, Equations \eqref{Stima1} and \eqref{Bound3} imply that
\begin{equation}\label{Bound4}
\mathbb{P}[Y=x|Y_1,Y_2,Y_3\in A_1]\leq \left(\frac{1}{1-\epsilon_4}\right)^3\mathbb{P}[\tilde{Y}=x]\leq \left(\frac{1}{1-\epsilon_4}\right)^3 \frac{C_2\lambda}{(1-\epsilon_{5})}.
\end{equation}
Here the last inequality holds because, since $\lambda\leq 9/10$, $\epsilon_5<1/10$ and $|A_1|$ is integer, $|A_1|\geq 2$ and we can apply Proposition \ref{Uniform} to the distribution $\tilde{Y}$.

Recalling that $$\mathbb{P}[Y_1\in A_1]\geq |A_1| \lambda (1-\epsilon_{4})\geq (1-\epsilon_5)(1-\epsilon_4)>(1-\epsilon_4-\epsilon_5),$$
we have $\mathbb{P}[Y_1\not\in A_1]\leq \epsilon_{4}+\epsilon_5$.
Proceeding as in Equation \eqref{case1b}, it follows that $$\mathbb{P}[Y=x]\leq \mathbb{P}[Y=x|Y_1,Y_2,Y_3\in A_1]+3\mathbb{P}[Y=x|Y_1\not\in A_1]\mathbb{P}[Y_1\not\in A_1]\leq$$
\begin{equation}\label{case1c}
\left(\frac{1}{1-\epsilon_4}\right)^3 \frac{C_2\lambda}{(1-\epsilon_{5})}+3\lambda\mathbb{P}[Y_1\not\in A_1]\leq \left(\frac{1}{1-\epsilon_4}\right)^3 \frac{C_2\lambda}{(1-\epsilon_{5})}+3\lambda(\epsilon_{4}+\epsilon_5)=C_3 \lambda
\end{equation}
which concludes CASE 1.

CASE 2: $|A_1|< \frac{(1-\epsilon_{5})}{\lambda}.$
Here we prove that $\mathbb{P}[Y_1+Y_2=x]\leq C_3 \lambda.$
Indeed we have that
$$\mathbb{P}[Y_1+Y_2=x]=\sum_{y_1\in A_1}\mathbb{P}[Y_1=y_1]\mathbb{P}[Y_2=x-y_1]+\sum_{y_1\not\in A_1}\mathbb{P}[Y_1=y_1]\mathbb{P}[Y_2=x-y_1].$$
Here we note that
\begin{equation}
\mathbb{P}[Y_1=y_1]\leq \begin{cases}
\lambda \mbox{ if } y_1\in A_1;\\
\lambda(1-\epsilon_{4}) \mbox{ if } y_1\not\in A_1.
\end{cases}
\end{equation}
Therefore
$$\mathbb{P}[Y_1+Y_2=x]\leq \lambda \left(\sum_{y_1\in A_1}\mathbb{P}[Y_2=x-y_1]+ (1-\epsilon_{4})\sum_{y_1\not\in A_1}\mathbb{P}[Y_2=x-y_1]\right)=$$
\begin{equation}\label{case2c}
=\lambda - \lambda\epsilon_{4}\sum_{y_1\not\in A_1}\mathbb{P}[Y_2=x-y_1].
\end{equation}
Now we note that, named $A_2:=\{y_1\in \mathbb{Z}_p:\ y_1\not \in A_1\}$, $|A_2|\geq p-\frac{(1-\epsilon_{5})}{\lambda}$ and we have that
$$\mathbb{P}[Y_2\in A_2]=1- \mathbb{P}[Y_2\in A_1]\geq 1- \lambda|A_1|\geq \epsilon_5.$$
Also, given another set $B$ of the same cardinality, we have that
$$\mathbb{P}[Y_2\in B]\geq \mathbb{P}[Y_2\in A_2]\geq \epsilon_{5}.$$
Therefore, considering that $|\{z\in \mathbb{Z}_p:\ z=x-y_1,\ y_1\not\in A_1\}|=|A_2|$, Equation \eqref{case2c} can be written as
$$\mathbb{P}[Y_1+Y_2=x]\leq\lambda-\lambda\epsilon_4\mathbb{P}[Y_2\in\{z\in \mathbb{Z}_p:\ z=x-y_1,\ y_1\not\in A_1\}]\leq \lambda - \lambda\epsilon_{4}\epsilon_{5}=C_3\lambda.$$
Now it is enough to note that, since $Y_1,Y_2$ and $Y_3$ are independent,
$$\mathbb{P}[Y_1+Y_2+Y_3=x]=\sum_{y_3\in \mathbb{Z}_p}\mathbb{P}[Y_3=y_3]\mathbb{P}[Y_1+Y_2=x-y_3]\leq \sum_{y_3\in \mathbb{Z}_p}\mathbb{P}[Y_3=y_3]C_3\lambda\leq C_3\lambda$$
which concludes CASE 2.

The thesis follows showing, with Mathematica, that we can choose $C_3< 1-2.27\cdot 10^{-12}.$
\endproof
\begin{rem}
In Theorem \ref{General3}, we may also weaken the hypothesis that $\lambda<9/10$ and assume that $\lambda<1$. With the same proof, choosing $\epsilon_5<1-\lambda$, we find the existence of constant $C_3(\lambda)<1$ also for values of $\lambda$ that are close to one. On the other hand, this constant cannot be made explicit and it depends on $\lambda$. So, in order to obtain an absolute constant, we rather prefer to assume $\lambda$ smaller than a given value (i.e. $9/10$).
\end{rem}
As a consequence, we can state the following result which is analogous, for the case of $\mathbb{Z}_p$, to that of the previous section.
\begin{thm}\label{Generalell}
Let us consider $\lambda\leq \frac{9}{10}$, $\epsilon>0$, $p>\frac{2}{\lambda\epsilon}$ and let $Y=Y_1+\dots+Y_{\ell}$ where $Y_i$ are i.i.d. whose distributions are upper-bounded by $\lambda$.
Then, if $\ell$ is sufficiently large with respect to $\epsilon$, we have that
$$\left(\max_{x\in \mathbb{Z}_p} \mathbb{P}[Y=x]\right)\leq \epsilon \lambda.$$
In particular, if $\epsilon=\left(\frac{3}{\ell_0}\right)^\nu$ where $\nu:=-\log_3(C_3)>0$\footnote{Here the best approximation we have for this constant is $\nu> 2.06\cdot10^{-12}.$}, we can take $\ell\geq \ell_0$.
\end{thm}
\proof
First of all, we prove, by induction on $k$, that, assuming $p>\frac{2}{C_3^{k-1}\lambda}$,
\begin{equation}\label{induction}\left(\max_{x\in \mathbb{Z}} \mathbb{P}[Y_1+Y_2+\dots+Y_{3^k}=x]\right)\leq C_3^k\lambda.\end{equation}
BASE CASE. The case $k=1$ follows from Theorem \ref{General3}.

INDUCTIVE STEP. We assume $p>\frac{2}{C_3^{k}\lambda}$ and that Equation \eqref{induction} is true for $k$ and we prove it for $k+1$.
At this purpose we set $\tilde{Y}_1=Y_1+Y_2+\dots+Y_{3^k}$ and we note that, set $\tilde{\lambda}=C_3^k\lambda$ for the inductive hypothesis,
$$\mathbb{P}[\tilde{Y}_1=x]\leq \tilde{\lambda}.$$
Therefore, because of Theorem \ref{General3}, we have that, if $\tilde{Y}_2$, $\tilde{Y}_3$, are three identical copies of $\tilde{Y_1}$ and set $Y= \tilde{Y}_1 +\tilde{Y}_2+\tilde{Y}_3$,
$$\mathbb{P}[Y=x]\leq C_3\tilde{\lambda}=C_3^{k+1}\lambda.$$
The inductive claim follows since $Y= \tilde{Y}_1 +\tilde{Y}_2+\tilde{Y}_3=Y_1+Y_2+\dots+Y_{3^{k+1}}$.

Now we consider $\ell \geq 3^k$. A bound for this case follows since we have that
$$\mathbb{P}[Y=x]=\sum_{y_3}\mathbb{P}[Y_{3^k+1}+Y_{3^k+2}+\dots+Y_{\ell}=y]\mathbb{P}[Y_1+Y_2+\dots +Y_{3^k}=x-y]\leq$$
$$ \sum_{y}\mathbb{P}[Y_{3^k+1}+Y_{3^k+2}+\dots+Y_{\ell}=y] C_3^k\lambda\leq C_3^k\lambda.$$

We have proved the claimed bound when $\epsilon=C_3^k$, $\ell\geq \ell_0=3^k$, and $p>\frac{2}{C_3^{k-1}\lambda}$. Here $\ell_0$ is such that $k=\lfloor \log_3(\ell_0)\rfloor$.
Therefore, set $\nu=-\log_3(C_3)$, we have
$$C_3^k=C_3^{\lfloor \log_3(\ell_0)\rfloor}=(3^{\log_3(C_3)})^{\lfloor \log_3(\ell_0)\rfloor}=(3^{\lfloor \log_3(\ell_0)\rfloor})^{-\nu}<(3^{\log_3(\ell_0)-1})^{-\nu}=\left(\frac{3}{\ell_0}\right)^\nu.$$
Also
$$p>\frac{2}{C_3^{k-1}\lambda}=\frac{2C_3}{(3^{\lfloor \log_3(\ell_0)\rfloor})^{-\nu}\lambda}=\frac{2C_3(3^{\lfloor \log_3(\ell_0)\rfloor})^{\nu}}{\lambda}$$ holds when
$$p>\frac{2C_3(\ell_0)^{\nu}}{\lambda}.$$
This implies that, for $\epsilon=\left(\frac{3}{\ell_0}\right)^\nu$, we can take any $\ell\geq\ell_0$ and
$$p>\frac{2C_3(\ell_0)^{\nu}}{\lambda}=\frac{2C_33^{\nu}}{\epsilon\lambda}=\frac{2}{\epsilon\lambda}.$$
\endproof
\section{A Littlewood-Offord Kind of Problem}
Now we will use the results of the previous sections to prove the following Littlewood-Offord kinds of results.
\begin{thm}
Let $A=\{v_1,v_2,\dots,v_n\}$ be distinct integers, and let us choose, uniformly at random, a subset $X$ of $A$ whose size is $\ell$. Then, however we take $\epsilon$, if $\ell$, $n$ and $n-\ell$ are sufficiently large, we have that
$$\left(\max_{x\in \mathbb{Z}} \mathbb{P}[\sum_{v\in X} v=x]\right)\leq\frac{\epsilon}{n}.$$
\end{thm}
\proof
We set $Y=Y_1+\dots+Y_{\ell}$ where $Y_i$ are independent and uniform distribution over $A$. Then
$$ \mathbb{P}[Y=x]=\mathbb{P}[Y=x|Y_i\not=Y_j,\ \forall i\not=j]\mathbb{P}[Y_i\not=Y_j,\ \forall i\not=j]+ $$
$$\mathbb{P}[Y=x|Y_i=Y_j \mbox{ for some } i\not=j]\mathbb{P}[Y_i=Y_j \mbox{ for some } i\not=j].$$
Here we note that, chosen, uniformly at random $X\subseteq A$ of size $\ell$, we have that
$$\mathbb{P}[\sum_{v\in X} v=x]=\mathbb{P}[Y=x|Y_i\not=Y_j,\ \forall i\not=j].$$
Hence we obtain the following relation, between $Y$ and the probability that $X$ sum to $x$
\begin{equation}\label{relation1} \mathbb{P}[Y=x]\geq \mathbb{P}[\sum_{v\in X} v=x]\mathbb{P}[Y_i\not=Y_j,\ \forall i\not=j].\end{equation}
Here we note that, if $n$ is sufficiently large with respect to $\ell$, we have $\mathbb{P}[Y_i\not=Y_j,\ \forall i\not=j]>1/2$. Therefore, for $n$ sufficiently large (with respect to $\ell$), we obtain from Equation \eqref{relation1} that
\begin{equation}\label{relation2} \mathbb{P}[Y=x]\geq \frac{1}{2}\mathbb{P}[\sum_{v\in X} v=x].\end{equation}
Because of Theorem \ref{Z,ell,epsilon}, however we consider $\epsilon$, if $\ell$ is sufficiently large (with respect to $\epsilon'=\epsilon/2$) and if $n$ is sufficiently large (with respect to $\ell$), we obtain that Equation \eqref{relation2} can be written as
\begin{equation}\label{relation3} \frac{\epsilon}{n}=\frac{2\epsilon'}{n}\geq2 \mathbb{P}[Y=x]\geq \mathbb{P}[\sum_{v\in X} v=x].\end{equation}
In the following part of this proof, we want to improve the result of Equation \eqref{relation3}. More precisely, we want to obtain it for any $n$ bigger than a value that depends only on $\epsilon$ (and does not depend on $\ell$).

Given $\epsilon$, note that Equation \eqref{relation3} holds in particular if we consider the smallest $\ell$ for which Theorem \ref{Z,ell,epsilon} holds, say $\ell_0$, and for any $n$ large enough with respect to $\ell_0$, say $n>n_0$. Let now consider $n>2n_0$ and let $\frac{n}{2}\geq \ell\geq\ell_0$. Here, chosen, uniformly at random $X\subseteq A$ of size $\ell$, $X_1\subseteq A$ of size $\ell-\ell_0$ and $X_2\subseteq (A\setminus X_1)$ of size $\ell_0$, we have that
$$\mathbb{P}[\sum_{v\in X} v=x]=\sum_{y\in \mathbb{Z}}\mathbb{P}[\sum_{v\in X_1} v=y]\mathbb{P}[\sum_{w\in X_2} w=x-y].$$
Since $|X_2|=\ell_0$ and
$$|A\setminus X_1|=n-(\ell-\ell_0)\geq n-\ell\geq n/2\geq n_0,$$
we have that $\mathbb{P}[\sum_{w\in X_2} w=x-y]\leq \frac{\epsilon}{n}$.
It follows that
\begin{equation}\label{relation4}\mathbb{P}[\sum_{v\in X} v=x]\leq \sum_{y\in \mathbb{Z}}\mathbb{P}[\sum_{v\in X_1} v=y]\frac{\epsilon}{n}=\frac{\epsilon}{n}.\end{equation}
This implies the thesis in case $\ell\leq n/2$. Note that this assumption is not a restriction since, however, we consider $X\subseteq A$ of size $\ell>n/2$, its complement $X'=A\setminus X$ has size $n-\ell\leq n/2$ and, set $x'=(\sum_{v\in A} v) -x$, $$\mathbb{P}[\sum_{v\in X} v=x]=\mathbb{P}[\sum_{w\in X'} w=x'].$$
\endproof
With essentially the same proof and by using Theorem \ref{Zk,ell,epsilon} instead of Theorem \ref{Z,ell,epsilon}, we obtain that:
\begin{thm}\label{LittlewoodOffordZk}
Let $A=\{v_1,v_2,\dots,v_n\}$ be distinct elements of $\mathbb{Z}_k$ where the prime factors of $k$ are larger than $n!$, and let us choose, uniformly at random, a subset $X$ of $A$ whose size is $\ell$. Then, however we take $\epsilon$, if $\ell$, $n$ and $n-\ell$ are sufficiently large, we have that
$$\left(\max_{x\in \mathbb{Z}} \mathbb{P}[\sum_{v\in X} v=x]\right)\leq\frac{\epsilon}{n}.$$
\end{thm}
And, using Theorem \ref{Generalell}:
\begin{thm}\label{LittlewoodOffordZp}
Let $A=\{v_1,v_2,\dots,v_n\}$ be distinct elements of $\mathbb{Z}_p$ and let us choose, uniformly at random, a subset $X$ of $A$ whose size is $\ell$. Then, however we take $\epsilon$, if $\ell$, $n$ and $n-\ell$ are sufficiently large and if $\frac{4n}{\epsilon}<p$, we have that
$$\left(\max_{x\in \mathbb{Z}_p} \mathbb{P}[\sum_{v\in X} v=x]\right)\leq\frac{\epsilon}{n}.$$
In particular, if $\epsilon=2\left(\frac{3}{\ell_0}\right)^\nu$ where $\nu=-\log_3(C_3)>0$, we can take $n-\ell$ and $\ell>\ell_0$ and $p>n^2$.
\end{thm}
Here the hypothesis $p>n^2$ is sufficient since we can assume $n/2\geq \ell_0$ and hence
$$n^2>n2\ell_0>n2\left(\frac{\ell_0}{3}\right)^{\nu}=\frac{4n}{\epsilon}.$$
\section{An application: $\Gamma$-sequenceability}
This section aims to consider a variation of the set-sequenceability problem. Here we introduce the following concept of $\Gamma$-sequenceability.
\begin{defn}
Let $\Gamma$ be a graph whose vertex set is $\{1,2,\dots, n\}$.
Then, a subset $A$ of size $n$ of a group $G$ is $\Gamma$-{\em sequenceable} if there is an ordering $(v_1, \ldots, v_n)$ of its elements such that the partial sums~$(s_0, s_1, \ldots, s_n)$, given by $s_0 = 0$ and $s_i = \sum_{j=1}^i v_j$ for $1 \leq i \leq n$, are different whenever $\{i,j\}\in E(\Gamma)$.

An ordering of $A$ for which $s_i\not=s_j$ whenever $\{i,j\}\in E(\Gamma)$, is said to be a $\Gamma$-sequencing.
\end{defn}
This concept generalizes that of $t$-{\em weak sequenceability}. Indeed, if $$E(\Gamma)=\{\{i,j\}: i,j\in [1,n], i\not=j,\ |j-i|\leq t\},$$
the concept of $\Gamma$-sequenceability reduces to that of $t$-weak sequenceability.

We recall that, in \cite{CD1}, the authors introduced a hybrid Ramsey theoretical and Probabilistic approach to proving that, when $n$ is large enough (with respect to $t$), a subset of a generic group $G$ is $t$-weak sequenceable. In this section, we will revisit their approach in the case $G=\mathbb{Z}_p$. In particular, using the improved bounds of the previous sections, we will prove that, given a graph $\Gamma$ whose maximum degree is at most $d$, then, any subset of a $\mathbb{Z}_p$ is $\Gamma$-sequenceable provided that $n$ is large enough with respect to $d$.

First of all, we recall the following proposition (see Proposition 3.2 of \cite{CD1}):
\begin{prop}\label{Ramsey}
Let $A$ be a subset of size $n$ of a (not necessarily abelian) group $G$ and let $t$ and $m$ be positive integers. Then, there exists a constant $n_{d,m}$ such that, if $n>n_{t,m}$, $A$ contains a subset $T$ whose size is at least $m$ that does not admit zero-sum subsets (with respect to any ordering) of size $g\leq t$.
\end{prop}
Also, with essentially the same inductive proof of Proposition 2.2 of \cite{CD}, we can prove that:
\begin{prop}\label{fix3}
Let $A\subseteq G\setminus\{0\}$ be a set of size $n$, let $\Gamma$ be a graph whose vertex-set is $[1,n]$ and whose maximum degree is at most $d$ and let $h$ be a positive integer such that $h\leq n-d-1$. Then there is an ordering of $h$-elements of $A$ that we denote with $(v_1,\dots,v_h)$, such that
\begin{itemize}
\item[(a)] for any ${0\leq i<j\leq h}$
$$s_i\not=s_j \mbox{ whenever } \{i,j\}\in E(\Gamma);$$
\item[(b)] for any $0\leq i \leq h$
$$s_i\not=s_n=\sum_{v\in A}v.$$
\end{itemize}\end{prop}

Then, using a probabilistic approach, we are able to prove the following theorem.
\begin{thm}
Let $\Gamma$ be a graph whose vertex-set is $[1,n]$ and whose maximum degree is at most $d$.
Then any subset of size~$n$ of~$\mathbb{Z}_p\setminus\{0\}$ is $\Gamma$-sequenceable whenever $n$ is large enough with respect to $d$ and $p>n^2$.
\end{thm}
\proof
Let $B$ be a generic subset of size $\bar{n}$ of $\mathbb{Z}_p\setminus \{0\}$.
Because of Theorem \ref{LittlewoodOffordZp}, we can fix $t>d$ so that, if $\bar{n}$ is large enough (i.e. $\bar{n}>n_0$), and if $p>\bar{n}^2$, the probability that a randomly chosen subset $X$ of size $\ell\in[t,\bar{n}-t]$ of $B$ sums to a given value $x$ is at most $\frac{1}{3d\bar{n}}$.

Now we fix a subset $A$ of size $n$ of $\mathbb{Z}_p\setminus \{0\}$.
Then, according to Proposition \ref{Ramsey}, however we choose $m$, if we assume that $n$ is large enough, there exists a subset $T$ of $A$ whose size is $m$ and that does not admit zero-sum subset of size $g\leq t$. Let us set $U=A\setminus T$ and let $n'=n-m$ be the size of $U$. According to Proposition \ref{fix3}, given $h= n'-t=n-m-t$, we can order $h$ elements of $U$, namely $(v_1,\dots,v_h)$, in such a way that
\begin{itemize}
\item[(a)] for any ${0\leq i<j\leq h}$
$$s_i\not=s_j \mbox{ whenever } \{i,j\}\in E(\Gamma);$$
\item[(b)] for any $0\leq i \leq h$
$$s_i\not=s_n=\sum_{v\in A}v.$$
\end{itemize}
We denote by $U'$ the set $\{v_1,\dots,v_h\}$.
Now it suffices to order the remaining $t$ elements of $U$ and the $m$ elements of $T$. Here we choose, uniformly at random an ordering $v_1,v_2,\dots,v_h,$ $z_{h+1},\dots, z_{n}$ that extend the one of $U'$. Let $X$ be the random variable that represents the number of ordered pairs $(i, j)$ such that $s_i = s_j$ with $0\leq i<j \leq n$ and $\{i,j\}\in E(\Gamma)$. We evaluate the expected value of $X$.

Following \cite{CD1}, because of the linearity of the expectation, we have that:
\begin{equation}\label{E1}\mathbb{E}(X)=\sum_{\substack{0\leq i<j \leq h\\ \{i,j\}\in E(\Gamma)}} \mathbb{P}(s_i=s_j)+\sum_{\substack{0\leq i\leq h<j\leq h+t\\ \{i,j\}\in E(\Gamma)}} \mathbb{P}(s_i=s_j)+\sum_{\substack{\max(h+t,i)<j\leq n\\ \{i,j\}\in E(\Gamma)}} \mathbb{P}(s_i=s_j).\end{equation}
Due to the choice of $U'$, when $0\leq i<j \leq h$, the probability $\mathbb{P}(s_i=s_j)=0$. Assuming that $0\leq i\leq h<j\leq h+t$,
since the equation
$$s_i=v_1+v_2+\dots+v_i =v_1+v_2+\dots+v_{h}+z_{h+1}+\dots+z_{j-1}+x=s_j$$ admits only one solution, $\mathbb{P}(s_i=s_j)$ is smaller or equal to $\frac{1}{|A\setminus (U'\cup\{z_{h+1},\dots,z_{j-1}\})|}\leq 1/m$ because $j$ is at most $h+t$. It follows that Equation \eqref{E1} can be rewritten as
\begin{equation}\label{E2}\mathbb{E}(X)\leq 0+\frac{t^2}{m}+\sum_{\substack{\max(h+t,i)<j\leq n\\ \{i,j\}\in E(\Gamma)}} \mathbb{P}(s_i=s_j).\end{equation}
Now, contrary to what was done in \cite{CDOR}, we need to divide the set of pairs $\mathcal{P}:=\{(i,j):\ \{i,j\}\in E(\Gamma): \max(h+t,i)<j\leq n\}$ in three parts:
\begin{itemize}
\item[(1)] $\mathcal{P}_1:=\{(i,j)\in \mathcal{P}: i<j\mbox{ and } j-i\leq t\}$;
\item[(2)] $\mathcal{P}_2:=\{(i,j)\in \mathcal{P}: i<j\mbox{ and } n-t>j-i>t\}$;
\item[(3)] $\mathcal{P}_3:=\{(i,j)\in \mathcal{P}: i<j\mbox{ and } j-i\geq n-t\}$.
\end{itemize}
It is easy to see that, for $(i,j)\in \mathcal{P}_1$, $s_i=s_j$ whenever
$$z_{i+1}+z_{i+2}+\dots+z_j=0.$$
Reasoning as before we have that, if at least one $z_{i+1},z_{i+2},\dots,z_j$ is not in $T$, the probability that $z_{i+1}+z_{i+2}+\dots+z_j=0$ is at most $\frac{1}{m}$. Indeed $s_i\not=s_{i+1}$ and so we may assume $j\geq i+2$. Assuming $z_{i+1}\not\in T$, we have that
$$z_{i+1}+z_{i+2}+\dots+z_{j-1}+x=0$$
has at most one solution where $x$ varies among at least $m$ values.

Therefore we have that:
$$\mathbb{P}(z_{i+1}+z_{i+2}+\dots+z_j=0)\leq \frac{1}{m}\mathbb{P}(\{z_{i+1},z_{i+2},\dots,z_j\}\not\subseteq T)+$$ $$\mathbb{P}(z_{i+1},z_{i+2},\dots,z_j\in T)\mathbb{P}(z_{i+1}+z_{i+2}+\dots+z_j=0|z_{i+1},z_{i+2},\dots,z_j\in T)$$
and, since $z_{i+1}+z_{i+2}+\dots+z_j\not=0$ whenever $z_{i+1},z_{i+2},\dots,z_j\in T$,
$$\mathbb{P}(z_{i+1}+z_{i+2}+\dots+z_j=0)\leq \frac{1}{m}\mathbb{P}(\{z_{i+1},z_{i+2},\dots,z_j\}\not\subseteq T).$$
Since $z_{i+1},z_{i+2},\dots,z_j$ are at most $t$ elements randomly chosen in $S\setminus U'$, the probability that all of them are contained in $T$ is at least $$\frac{{|T|\choose j-i}}{{|S\setminus U'|\choose j-i}}\geq\frac{{|T|\choose t}}{{|S\setminus U'|\choose t}} = \frac{{m \choose t}}{{m+t\choose t}}\geq \frac{(m-(t-1))^{t}}{(m+t)^{t}}.$$
It follows that, if $(i,j)\in \mathcal{P}_1$,
$$\mathbb{P}(z_{i+1}+z_{i+2}+\dots+z_j=0)\leq \frac{1}{m}\left(1-\frac{(m-(t-1))^{t}}{(m+t)^{t}}\right).$$
Let now consider $(i,j)\in \mathcal{P}_2$. Here we have that
$$s_i=\begin{cases}v_1+v_2+\dots+v_i \mbox{ if }i<h;\\
v_1+v_2+\dots+v_{h}+z_{h+1}+\dots+z_{i-1}+z_i\mbox{ otherwise}.
\end{cases}
$$
Hence $s_i=s_j$ implies that
$$\begin{cases}z_{h+1}+\dots+z_{j-1}+z_j=-(v_{i+1}+\dots +v_{h}) \mbox{ if }i<h;\\
z_{i+1}+\dots+z_{j-1}+z_j=0\mbox{ otherwise}
\end{cases}
$$
Hence, since both $j-i$ and $j-h$ are at least $t$, given $(i,j)\in \mathcal{P}_2$, and assuming $|A\setminus U'|=m+t$ is large enough (i.e. $m+t>n_0$), due to Theorem \ref{LittlewoodOffordZp} and our assumption on $t$, we have that
\begin{equation}\label{p2}\mathbb{P}[s_i=s_j] \leq \frac{1}{3d(m+t)}\leq \frac{1}{3dm}.\end{equation}
Finally, we consider $(i,j)\in \mathcal{P}_3$. Here, if $n$ is large enough, we have that $n-m-t=h>t$ and hence
$$n-i\geq j-i\geq n-t> n-h.$$
It follows that we may assume $i< h$. Also, since $s_i\not=s_n$ for any $i< h$, it suffices to consider $h<j<n$. Here $s_i=s_j$ implies that
$$s_i=v_1+\dots+v_i=v_1+\dots+v_h+z_{h+1}+\dots+z_j=s_j$$
and, hence
$$v_{i+1}+\dots+v_h+z_{h+1}+\dots+ z_j=0.$$
Since we are assuming $j<n$, we also have that
$$z_{j+1}+\dots+z_n+v_{1}+\dots+v_{i}=s_n$$
and, considering that
$$z_{j+1}+\dots+x+v_{1}+\dots+v_{i}=s_n$$
has exactly one solution, and $n-j<n-(j-i)\leq t$, given $(i,j)\in \mathcal{P}_3$,
\begin{equation}\label{p2}\mathbb{P}[s_i=s_j] \leq \frac{1}{m}.\end{equation}

Hence Equation \eqref{E2} becomes
\begin{equation}\label{E3}\mathbb{E}(X)\leq \frac{t^2}{m}+\sum_{(i,j)\in \mathcal{P}_1} \frac{1}{m} \left(1-\frac{(m-(t-1))^{t}}{(m+t)^{t}}\right)+\sum_{(i,j)\in \mathcal{P}_2}\frac{1}{3dm}+\sum_{(i,j)\in \mathcal{P}_3}\frac{1}{m}.\end{equation}

Now we note that there are less than $(m+t)t$ pairs $(i,j)\in \mathcal{P}_1$. Also, there are at most $(n-h-t)=m$ values of $j$ such that $h+t<j\leq n$ and hence there are at most $md$ pairs $(i,j)\in \mathcal{P}_2$. Given $(i,j) \in \mathcal{P}_3$, $j$ must be larger than $n-t$ and so there are at most $td$ pairs $(i,j)\in \mathcal{P}_3$.
Therefore, from Equation \eqref{E3} we obtain
\begin{equation}\label{E3}\mathbb{E}(X)\leq \frac{t^2}{m}+\frac{(m+t)t}{m} \left(1-\frac{(m-(t-1))^{t}}{(m+t)^{t}}\right)+\frac{1}{3}+\frac{td}{m}.\end{equation}
Here we note that the right-hand side of the last inequality goes to $1/3$ as $m$ goes to infinite. Therefore for $m$ large enough, or more precisely for $m\geq m_0$ for a suitable $m_0\in \mathbb{N}$, we have that $\mathbb{E}(X)<1$. This means that, if $n$ is large enough (i.e. if $n>n_{t,m_0}$) there exists an ordering on which $X=0$ and hence there exists a $\Gamma$-sequencing of $A$.
\endproof
Note that, a similar result can be obtained also assuming that $A$ is a subset of $\mathbb{Z}_k$ and the prime factors of $k$ are large enough. Here we need to use Theorem \ref{LittlewoodOffordZk} instead of Theorem \ref{LittlewoodOffordZp}, and we get that:
\begin{thm}
Let $\Gamma$ be a graph whose vertex-set is $[1,n]$ and whose maximum degree is at most $d$.
Then any subset of size~$n$ of~$\mathbb{Z}_k\setminus\{0\}$ is $\Gamma$-sequenceable whenever $n$ is large enough with respect to $d$ and the prime factors of $k$ are larger than $n!$.
\end{thm}

\section*{Acknowledgements}
The author would like to thank Stefano Della Fiore for our useful discussions on this topic. The author was partially supported by INdAM--GNSAGA.

\end{document}